\documentclass[letterpaper, 10 pt, conference]{ieeeconf}  

\IEEEoverridecommandlockouts                              

\usepackage{graphics} 
\usepackage{epsfig} 
\usepackage{mathptmx} 
\usepackage{times} 
\usepackage{amsmath} 
\usepackage{amssymb}  
\usepackage{mdwlist}
\usepackage{afterpage}
\usepackage{subfigure}
\usepackage{color,soul}
\newtheorem{theorem}{Theorem}[section]

\newtheorem{proposition}[theorem]{Proposition}

\newenvironment{remark}[1][Remark]{\begin{trivlist}
\item[\hskip\labelsep{\indent\itshape #1:}]}{\end{trivlist}}

\title{\LARGE \bf
Shaping up crowd of agents through controlling their statistical moments
}

\author{Yuecheng Yang$^{\text{ a,b}}$, Dimos V. Dimarogonas$^{\text{ a,c}}$ and Xiaoming Hu$^{\text{ a,b}}$%
\thanks{\hspace{-3mm}$^{\text{ a}}$ ACCESS Linnaeus Centre, KTH - Royal institute of Technology\newline
$^{\text{ b}}$ Automatic Control Lab, School of Electrical Engineering, KTH - Royal institute of Technology\newline
$^{\text{ c}}$ Centre for Autonomous Systems, KTH - Royal institute of Technology\newline
$^{\text{ d}}$ Optimization and Systems Theory, Department of Mathematics, KTH - Royal institute of Technology\newline
This work is supported by SSF, the Swedish Research Council}%
}

\begin{document}

\maketitle
\thispagestyle{empty}
\pagestyle{empty}

\begin{abstract}
In a crowd model based on leader-follower interactions, where positions of the leaders are viewed as the control input, up-to-date solutions rely on knowledge of the agents' coordinates. In practice, it is more realistic to exploit knowledge of statistical properties of the group of agents, rather than their exact positions. In order to shape the crowd, we study thus the problem of controlling the moments instead, since it is well known that shape can be determined by moments. An optimal control for the moments tracking problem is obtained by solving a modified Hamilton-Jacobi-Bellman (HJB) equation, which only uses the moments and leaders' states as feedback. The optimal solution can be solved fast enough for on-line implementations.
\end{abstract}

\section{Introduction}
With continuous urbanization of the global population, researchers are getting more and more aware that it is important to have a better understanding of crowd behavior under certain circumstances. Experts from different fields, including sociologists, psychologists, ecologists, physicists, mathematicians and computer scientists use their different perspectives to model and analyze the crowd. One typical and important application of the crowd behavior study is the evacuation problem in emergency situations. Although the human behavior can be very complicated in these situations, researchers have made lots of effort to model, analyze and simulate the crowd with different approaches for the purpose of minimizing the total societal loss.

From the pure social psychology point of view, there are models such as the theory of planned behavior form I. Ajzen \cite{Ajzen1991}. J.D. Sime linked the psychology part and the engineering part of the crowd behavior problem in \cite{Sime1995}. From the engineering aspect, one commonly used macroscopic approach is to consider the whole crowd as one entity that is described by a density function and use tools such as fluid mechanics and partial differential equations to conduct analysis. This continuum setting of the crowd has been used, e.g., in \cite{continuum} and \cite{flow}. The disadvantage of this type of approach relies on the fact that a density approximation of the human crowd is not adequate when the density is low. On the other hand, low computational cost is obviously the advantage of the approach, especially when the scale of the crowd is large enough.

Another method for crowd analysis is by using the tool of multi-agent systems theory that has rapidly developed in the past decades. Compared to the continuum setting, multi-agent system analysis can be seen as a microscopic model that focuses on individual behaviors. Part of the idea came from animal flocking observation and modeling such as Reynolds model in \cite{reynolds} in 1987. This kind of approach was used, e.g., in \cite{simulating}, \cite{Helbing2009}. One of the most widely used models, the social force model in \cite{Social force}, also uses this setting. With the theory development in multi-agent consensus problems, researchers do strict analysis for linear multi-agent systems, e.g., in \cite{Olfati-Saber2003}. However, many models for human crowd are highly nonlinear, which makes the analysis much harder due to the lack of theoretical support for nonlinear systems. The computation cost is also a big issue for numerical experiments. Nevertheless, the advantage is the accuracy of individual states if compared to macroscopic approaches.

Among the multi-agent crowd models, there is one type of model called leader-follower model that divides the crowd into two groups based on their roles. These models are very useful for example in the evacuation problem, where the rescue workers act as leaders and the general public can be considered as followers. Such models can be found for example in \cite{Morse2003}, \cite{Eren2005} and \cite{Olfati-Saber2006}.  A central issue for leader-follower models is controllability. In \cite{magnus2006}, the controllability of linear leader-follower models is discussed in detail. Unfortunately, the system becomes uncontrollable in most cases of realistic linear models. In \cite{cdc2013}, an optimal control approach is used for control design even though the system is not controllable for the evacuation problem. In leader-follower based models, the formation of the followers is in general very hard to control by only the leaders, since the system is often uncontrollable. Unfortunately, the possibility of shaping the crowd is seldom discussed in the literature since the shape is less restrictive than the formation and maybe more elaborate to control. In this paper, we will study the moments of the crowd, which have a strong connection to the shape as showed in \cite{Milanfar1995} and \cite{Golub2004}. By using dynamic programming techniques, we will attack the question of shaping the crowd.

In the method introduced later in this article, we will design a controller so that the shape of the crowd will track approximately some desired shape during the process. Although the asymptotic behavior of the crowd is important, the transient phase is even more critical for crowd control problems. Lyapunov stability analysis does not cover the transient behavior of the system in a straightforward manner, thus is not used in the paper. We instead introduce a new control design approach by using optimal control theory. Moreover, the leaders will need only the moments information together with their own states to calculate their movement as time evolves even though the optimal control problem has all the followers' positions and velocities as its state variables. This make the method introduced in the paper more practical since individual states are in general very hard to measure. The numerical method introduced is shown to be efficient enough for on-line implementation as well.

The outline of this paper is as follows: in Section II, the basic leader-follower model is stated and the moments of the crowd are defined. There is a short discussion about the moments tracking problem in general and the motivation why the optimal control approach is used in this paper. The optimal control problem is introduced and approximately solved in Section III by using model predictive control tools. In Section IV, some improvements of the method are made to handle numerical issues while several detailed experiments are carried out to test the capability and robustness of the method. A short conclusion can be found in Section V together with a brief outline for future work.

\section{Problem formulation}
\subsection{Leader-follower model}
We will study the behavior of two types of agents in a two dimensional space. The first type of agents are called {\it followers}. We assume that there are $N$ followers in total, and their motion follows some simple rules that mostly depend on the relative positions and velocities among them and on those to the other type of the agents - the {\it leaders}. The $M$ leaders, on the other hand, have better knowledge of the whole environment and can plan their motion accordingly. Based on the situation of the whole crowd, the leaders should behave in an optimal way to guide the followers to reach certain goals. The position of follower $i$ is denoted by a vector $z_i\in \mathbb R^2$, and its velocity denoted by $v_i\in \mathbb R^2$. The position of leader $j$ is denoted by $z_{l_j}\in\mathbb R^2$, and its velocity by $v_{l_j}\in\mathbb R^2$.

Suppose that without the leaders, the influence from follower $k$ to follower $i$ is modeled as $f(z_i,v_i,z_k,v_k)$. By assuming unit mass for each follower, we can write
\begin{equation}
\dot v_i=\sum_{k=1\atop{k\ne i}}^{N} f(z_i,v_i,z_k,v_k)-p v_i,
\end{equation}
where $p$ is a coefficient to model the damping effect such as resistance or physical limit of the agent. When follower $i$ senses any of the leaders, the influence by the leader is added to the model:
\begin{equation}
\dot v_i=\sum_{j=1}^{M}g(z_i,v_i,z_{l_j},v_{l_j})+\sum_{k=1\atop k\ne i}^{N} f(z_i,v_i,z_k,v_k)-p v_i.
\end{equation}
The functions $f$ and $g$ will be defined or described later. Note that here we sum up all the terms with the $M$ leaders, meaning that the leader-follower interaction is state-based. There is no predesigned network. If $g$ has no zeros in its domain, then the followers can always sense all the leaders and be influenced more or less by all of them.

For the leader $j$, we have
\begin{equation}
\dot v_{l_j}=u_j-p v_j,
\end{equation}
so it has the same damping coefficient as the followers and $u_j$ is the input to be designed in order to fulfill certain tasks. We denote $$z=\begin{bmatrix}z_{1}\\z_{2}\\\vdots\\z_{N}\end{bmatrix},v=\begin{bmatrix}v_{1}\\v_{2}\\\vdots\\v_{N}\end{bmatrix},z_l=\begin{bmatrix}z_{l_1}\\z_{l_2}\\\vdots\\z_{l_M}\end{bmatrix},v_l=\begin{bmatrix}v_{l_1}\\v_{l_2}\\\vdots\\v_{l_M}\end{bmatrix}.$$

If there is a cost function defined over time and a terminal cost to be minimized,  then we can formulate an optimal control problem:
\begin{equation}
\begin{array}{cl}
\displaystyle\min_{u}& \Phi(t_f,z(t_f),v(t_f),z_l(t_f),v_l(t_f))\\\vspace{1mm}
&~~~~~~~~~+\displaystyle\int_{t_0}^{t_f}\mathcal L(t,z(t),v(t),z_l(t),v_l(t),u(t))dt\\\vspace{1mm}
s.t.~ &\dot z_i=v_i,\\\vspace{1mm}
&\dot v_i=\sum_{j=1}^{M}g(z_i,v_i,z_{l_j},v_{l_j})+\sum_{k=1}^{N} f(z_i,v_i,z_k,v_k)-p v_i,\\\vspace{1mm}
&\dot z_{l_j}=v_{l_j},\\ \vspace{1mm}
&\dot z_{l_j}=u_j-pv_j,\\\vspace{1mm}
&u_j\in U,~ z(t_0),v(t_0), z_l(t_0),v_l(t_0)\text{ given, }
\end{array}
\end{equation}
for $i=1,\cdots,N$, and $j=1,\cdots,M$, where the set $U\subset \mathbb R^2$ is the feasible that will be specified later. The solution to this problem $u=\begin{bmatrix}u_{1}&u_{2}&\cdots&u_{M}\end{bmatrix}^T$will give a proper control input for the leaders.

\subsection{Definitions of the moments}
In order to shape the crowd, we need at least some statistic properties of it. In this article, we make use of moments of the agents' position. There are two ways to define moments:
\subsubsection{Definition 1}
The \emph{(raw) moment} of {\it order} $k$ is defined by
\begin{equation}
\label{moments1.1}
M_{ab}=\frac{1}{N}\sum_{i=1}^{N}z_{i,x}^a z_{i,y}^b,
\end{equation}
where $z_{i,x}$ and $z_{i,y}$ are the $x$, $y$ coordinate for $z_i, and $$a, b\in\mathbb N$ with $a+b=k$.

\subsubsection{Definition 2}
The \emph{centralized moment} is defined by
\begin{equation}
\label{moments2.1}
\bar M_{10}=\frac{1}{N}\sum_{i=1}^{N}z_{i,x},\bar M_{01}=\frac{1}{N}\sum_{i=1}^{N}z_{i,y},
\end{equation}
and
\begin{equation}
\label{moments2.2}
\bar M_{ab}=\frac{1}{N}\sum_{i=1}^{N}(z_{i,x}-\bar M_{10})^a(z_{i,y}-\bar M_{01})^b, \text{ for }k>1, \text{ and }a+b=k.
\end{equation}

It is not hard to show that one can calculate all the centralized moments by knowing all the raw moments and vise versa. If not particularly stated, we will only use the raw moments in the rest of the paper.

The moments have a strong relationship with the shape of the crowd. If the shape and the distribution in the shape is known, one can approximate the moments by
$$M_{ab}=\int_{S}x^ay^b\rho(x,y)dxdy,$$
where $S$ is the shape, and $\rho(x,y)$ is the density function of the distribution. On the other hand, if all the moments $M_{ab}$ with $a+b\leq 2n-3$ are given, then one can approximate the shape, assuming it is convex and a uniform distribution, by a $n$-polygon by using a modified method similar to the one introduced in \cite{Milanfar1995}.

\subsection{Moments tracking}


It is noticeable that the scale of the leader-follower model:
\begin{align}
\begin{cases}
 \dot z_i&=v_i,\notag\\
\dot v_i&=\sum_{j=1}^{M}g(z_i,v_i,z_{l_j},v_{l_j})+\sum_{k=1}^{N} f(z_i,v_i,z_k,v_k)-p v_i,\notag\\
\dot z_{l_j}&=v_{l_j},\notag\\
\dot z_{l_j}&=u_j-pv_j,
\end{cases}
\end{align}
for $i=1,\cdots,N$, and $j=1,\cdots,M$ becomes very large when the number of followers increases. Hence, individual control will be difficult even if the whole system is controllable, which is not always the case. This gives a motivation to look into the statistical properties of the crowd instead. Because of the strong connection between the moments and the shape, we can setup a moments tracking problem if we want to do shape tracking. Namely, from a desired shape ``signal" and a desired density distribution, we can calculate the corresponding desired moments. If moments of the crowd could track the desired moments, then the crowd should more or less follow the desired shape. The more moments we track, the better the performance should be.

Theoretically, if we can derive the evolution of the moments $M_{ab}$ with $a+b\leq m$ and their time derivatives by functions only using $M_{ab}$, $\frac{d}{dt}M_{ab}, a+b\leq m$ and $z_l,v_l,u$, then we will obtain a, probably nonlinear, system with $M_{ab}$ and $\frac{d}{dt}M_{ab}, a+b\leq m$ and $z_l,v_l$ as its states $u$ as its control. Furthermore, if this system is controllable, them for any continuous moments signal, we can get a minimum energy feedback control to track those moments. However, it is impossible in general to write down this type of system without using any information from $z$  and $v$ defined previously if no approximation is made. Even if for some special $f$ and $g$ functions, one can write down the system of moments, it is inaccessible for almost all cases that we have examined. This is still true even if one simplifies the leader-follower dynamics into a single integrator system. The analysis of some specific models are made in \cite{axiv} and is omitted here due to the space limitation.

Another way to deal with the tracking problem of an uncontrollable system is by setting a cost function for the system and using optimal control techniques to solve it. A very standard cost function is the quadratic errors between the real moments and the desired moments signals. Given a sequence of nonnegative scalars $\{c_{ab}\}$ with $a+b\leq m$, we can setup the following optimal control problem:
\begin{equation}
\label{OCP221}
\begin{array}{cl}\vspace{2mm}
\displaystyle \min_{u} & \displaystyle\sum_{a+b\leq m}\displaystyle\int_{t_0}^{t_f}c_{ab}(M_{ab}(t)-M_{ab}^d(t))^2dt\\\vspace{1mm}
s.t.  &\dot z_i=v_i,\\\vspace{1mm}
&\dot v_i=\sum_{j=1}^{M}g(z_i,v_i,z_{l_j},v_{l_j})+\sum_{k=1}^{N} f(z_i,v_i,z_k,v_k)-p v_i,\\\vspace{1mm}
&\dot z_{l_j}=v_{l_j},\\\vspace{1mm}
&\dot z_{l_j}=u_j-pv_j,\\\vspace{1mm}
 &\|u_j(t)\|\leq u_{max},\\\vspace{1mm}
&z(t_0), v(t_0), z_l(t_0), v_l(t_0) \text{ given, }
\end{array}
\end{equation}
where $M_{ab}^d(t)$ is the given signal for $M_{ab}(t)$ to track.

The analytic solution is possible to find only if $f$ and $g$ are simple functions such as linear functions. One can use Pontryagin's minimum principle (PMP) to solve such type of problem similar to the method used in \cite{cdc2013}. However, these two functions are usually nonlinear in practical models, which implies that PMP is to hard to solve. For the dynamic programming approach, if we can find a cost to go function $J$ that satisfies the Hamilton-Jacobi-Bellman equation, then the optimal control can be calculated by partial derivatives of $J$. The HJB equation of (\ref{OCP221}) can be written as:
\begin{align}
-\frac{\partial J}{\partial t}=&\min_{\|u\|\leq u_{max}} \{ \sum_{a+b\leq m}c_{ab}(M_{ab}(t)-M_{ab}^d(t))^2+ \frac{\partial J}{\partial z}^T v\notag\\
+&\sum_{i=1}^N\frac{\partial J}{\partial v_i}(\sum_{j=1}^{M}g(z_i,v_i,z_{l_j},v_{l_j})+\sum_{k=1}^{N} f(z_i,v_i,z_k,v_k)-p v_i)\notag\\
+&  \frac{\partial J}{\partial z_l}^T v_l+\frac{\partial J}{\partial v_l}^T (u-p v_l)\},
\label{HJB1}
\end{align}
with the boundary condition \begin{equation}J(t_f,z,v,z_l,v_l)=0.\end{equation}
The minimum of the right-hand side in equation (\ref{HJB1}) is reached when $u_j^*=-u_{max}\lambda_j$ where $\lambda_j$ is the unit vector of $\frac{\partial J}{\partial v_{l_j}}$. If we plug this back to (\ref{HJB1}), we get
\begin{align}
-\frac{\partial J}{\partial t}=& \sum_{a+b\leq m} c_{ab}(M_{ab}(t)-M_{ab}^d(t))^2+ \frac{\partial J}{\partial z}^T v\notag\\
+&\sum_{i=1}^N\frac{\partial J}{\partial v_i}(\sum_{j=1}^{M}g(z_i,v_i,z_{l_j},v_{l_j})+\sum_{k=1}^{N} f(z_i,v_i,z_k,v_k)-p v_i)\notag\\
+&  \frac{\partial J}{\partial z_l}^T v_l-\sum_{j=1}^{M}u_{max}\|\frac{\partial J}{\partial v_{l_j}}\|-p \frac{\partial J}{\partial v_l}^T v_l,
\label{HJB2}
\end{align}

Unfortunately, this partial differential equation is also very difficult to solve since the number of variables is proportional to the number of followers, which means the complex of solving the PDE becomes very high. Meanwhile, even if the HJB equation is sovable, one needs all the followers' states to calculate the optimal control, which is very difficult to implement in practice. Hence we will introduce a suboptimal control by using moments information as feedback while the complexity of the algorithm is low enough.

\section{Feedback control using moments information only}

In this section, we will study a special problem, which gives a reasonable model for crowd behavior. Let us make the following assumption for $g$:
\begin{itemize}
{\item The function $g$ is composed of two parts, a position consensus part and a velocity alignment part. There are two functions that give different weight to these two parts which only depend on the positions of the leader and the follower, i.e.,
$$g(z_i,v_i,z_{l_j},v_{l_j})=g_1(z_i,z_{l_j})(z_{l_j}-z_i)+g_2(z_i,z_{l_j})(v_{l_j}-v_i),$$
where $g_1$ and $g_2$ are real valued functions. Ideally, both $g_1$ and $g_2$ should be functions of the distance $d=\|z_i-z_{l_j}\|$ to make the model reasonable in practice. Moreover, $g_1$ should be relatively large when $d$ is large while $g_2$ should dominate when $d$ is close to zero. This is due to the fact that catching up the leader is more important when the distance is large while moving in the same direction makes more sense when the distance is already short enough.}
\end{itemize}

With these assumptions we can rewrite the HJB equation (\ref{HJB2}) as
\begin{align}
-\frac{\partial J}{\partial t}=& \mathcal L(t) + \frac{\partial J}{\partial z}^T v+\sum_{i=1}^N\frac{\partial J}{\partial v_i}\Big(\sum_{j=1}^{M}g_1(z_i,z_{l_j})(z_{l_j}-z_i)\notag\\
&+g_2(z_i,z_{l_j})(v_{l_j}-v_i)+f_i-p v_i\Big)+  \frac{\partial J}{\partial z_l}^T v_l\notag\\
&-\sum_{j=1}^{M}u_{max}\|\frac{\partial J}{\partial v_{l_j}}\|-p \frac{\partial J}{\partial v_l}^T v_l,
\label{HJB3}
\end{align}
where the notations $\mathcal L(t)=\sum_{a+b\leq m} c_{ab}(M_{ab}(t)-M_{ab}^d(t))^2$ and $f_i=\sum_{k=1}^{N}f(z_i,v_i,z_k,v_k)$ are introduced for simplicity.
\subsection{HJB simplification}

In order to avoid using all the followers' states $z$ and $v$, we want to derive $J$ as a function of $M_{ab}$ with $a+b\leq m$, $\dot M_{10}$, $\dot M_{01}$, $z_l$ and $v_l$ only, which means that the followers' positions only appear in the moments and their velocities only appear in the time derivative of the first order moments. In this section, we will show that this can be achieved with an approximated HJB equation.

\begin{proposition}
The Hamilton-Jacobi-Bellman equation (\ref{HJB3}) can be approximated in such a way so that the cost-to-go function $J$ has $t, M_{ab}$ with $a+b\leq m$, $\dot M_{10}$, $\dot M_{01}$, $z_l$ and $v_l$ as its variables.
\end{proposition}

 Assuming that there is a function $\bar J(t,M_{ab}, \dot M_{10}, \dot M_{01}, z_l, v_l)$ that satisfies (\ref{HJB3}), then we have the following equations because of the chain rule:
$$\frac{\partial \bar J}{\partial z_i}=\sum_{a+b\leq m}\frac{\partial \bar J}{\partial M_{ab}}\frac{\partial M_{ab}}{\partial z_i}=\frac{1}{N}\sum_{a+b\leq m}\frac{\partial \bar J}{\partial M_{ab}}\begin{bmatrix}az_{i,x}^{a-1}z_{i,y}^b\vspace{2mm} \\bz_{i,x}^{a}z_{i,y}^{b-1}\end{bmatrix}.$$
$$\frac{\partial \bar J}{\partial v_i}=\frac{\partial \bar J}{\partial \dot M_{10}}\frac{\partial \dot M_{10}}{\partial v_i}+\frac{\partial \bar J}{\partial \dot M_{01}}\frac{\partial \dot M_{01}}{\partial v_i}=\frac{1}{N}\begin{bmatrix}\frac{\partial\bar J}{\partial \dot M_{10}}\vspace{2mm} \\\frac{\partial J}{\partial \dot M_{01}}\end{bmatrix}.$$
Then we can write the HJB equation for $\bar J$ as
\begin{align}
-\frac{\partial \bar J}{\partial t}= &\mathcal L(t) +\frac{1}{N}\sum_{a+b\leq m}\frac{\partial \bar J}{\partial M_{ab}}\sum_{i=1}^N\begin{bmatrix}az_{i,x}^{a-1}z_{i,y}^b &bz_{i,x}^{a}z_{i,y}^{b-1}\end{bmatrix}v_i\notag\\
&+\frac{1}{N}\begin{bmatrix}\frac{\partial \bar J}{\partial \dot M_{10}}&\frac{\partial \bar J}{\partial \dot M_{01}}\end{bmatrix} \Big(  \sum_{j=1}^{M} \big(\sum_{i=1}^N\big(g_1(z_i,z_{l_j})(z_{l_j}-z_i)\notag\\
&+g_2(z_i,z_{l_j})(v_{l_j}-v_i)\big)\big)+\sum_{i=1}^Nf_i(z,v)-p\sum_{i=1}^N v_i\Big)\notag\\
&+\frac{\partial \bar J}{\partial z_l}^T v_l-p \frac{\partial \bar J}{\partial v_l}^T v_l-\sum_{j=1}^{M}u_{max}\|\frac{\partial \bar J}{\partial v_{l_j}}\|.
\label{HJB4}
\end{align}

Unfortunately, there are still terms on the right-hand side of the equation that contain $z$ and $v$, which implies that the assumption does not hold in general. We need to approximate those terms in order to get rid of $z$ and $v$. In order to do the approximation, we need to make some further assumptions:
\begin{itemize}
\item{$ f(z_i,v_i,z_k,v_k)=- f(z_k,v_k,z_i,v_i)$, meaning that the follower-follower interactions are symmetric. As a result, we get $\sum_{i=1}^N f_i(z,v)=0$.}
\item{The terms of the form $\sum_{i=1}^N h(z_i)v_i$ can be approximated by $(\sum_{i=1}^N h(z_i))\dot M_1 $ for any scalar function $h(\cdot)$, where $\dot M_1=\begin{bmatrix}\dot M_{10} & \dot M_{01}\end{bmatrix}^T$. This leads to the following approximation $$\sum_{i=1}^N z_{i,x}^az_{i,y}^b v_i\approx (\sum_{i=1}^N z_{i,x}^az_{i,y}^b )\dot M_1=NM_{ab}\dot M_1,$$ and $$\sum_{i=1}^{N}(g_2(z_i,z_{l_j})(v_{l_j}-v_i))\approx\Big(\sum_{i=1}^{N}(g_2(z_i,z_{l_j})\Big)(v_{l_j}-\dot M_1)$$ If we regard $z$ and $v$ as random variables, then this assumption is equivalent to saying that $z$ and $v$ are independent.}
\item{The functions $g_1(z_1,z_2)$ and $g_2(z_1,z_2)$ can be approximated by polynomials of $z_1$ with degree less than $m$ for a given $z_2$, i.e.,
\begin{align*}
&g_1(z_1,z_2)=\sum_{a=0}^{m-1}\sum_{b=0}^{m-a-1}\alpha_{ab}(z_2)z_{1x}^az_{1y}^b,\\
&g_2(z_1,z_2)=\sum_{a=0}^{m-1}\sum_{b=0}^{m-a-1}\beta_{ab}(z_2)z_{1x}^az_{1y}^b.
\end{align*}
This approximation can be usually achieved by Taylor expansion around certain points if the function $g_1$ and $g_2$ are ``good" enough (smooth and not very steep). Once we can use polynomials to approximate $g_1$ and $g_2$, we will have
\begin{align*}
&\sum_{i=1}^Ng_1(z_i,z_{l_j})(z_{l_j}-z_i)\\
\approx&\sum_{i=1}^N \sum_{a=0}^{m-1}\sum_{b=0}^{m-a-1}\alpha_{ab}(z_{l_j})z_{i,x}^az_{i,y}^b(z_{l_j}-z_i)\\
=&\sum_{a=0}^{m-1}\sum_{b=0}^{m-a-1} \alpha_{ab}(z_{l_j}) \sum_{i=1}^N z_{i,x}^az_{i,y}^b(z_{l_j}-z_i)\\
=&N\sum_{a=0}^{m-1}\sum_{b=0}^{m-a-1} \alpha_{ab}(z_{l_j})  \left(M_{ab}z_{l_j}-\begin{bmatrix}M_{a+1,b}\\M_{a,b+1}\end{bmatrix}\right)
\end{align*}
And similarly we have

$$\sum_{i=1}^N g_2(z_i,z_{l_j})= N\sum_{a=0}^{m-1}\sum_{b=0}^{m-a-1} \beta_{ab}(z_{l_j})M_{ab}.$$}
\end{itemize}
With the above assumptions the HJB equation (\ref{HJB4}) can be finally written as:
\begin{align}
-\frac{\partial \tilde{J}}{\partial t}&= \mathcal L(t) +\frac{\partial \tilde{J}}{\partial z_l}^T v_l-p \frac{\partial \tilde{J}}{\partial v_l}^T v_l-u_{max}\|\frac{\partial \tilde{J}}{\partial v_l}\|\notag\\
&+\sum_{a=0}^{m}\sum_{b=0}^{m-a}\frac{\partial \tilde{J}}{\partial M_{ab}}\begin{bmatrix}aM_{a-1,b} &bM_{a,b-1}\end{bmatrix}\dot M_1\notag\\
&+\begin{bmatrix}\frac{\partial \tilde{J}}{\partial \dot M_{10}}&\frac{\partial \tilde{J}}{\partial \dot M_{01}}\end{bmatrix} \Big(  \sum_{j=1}^{M} \big( \sum_{a=0}^{m-1}\sum_{b=0}^{m-a-1}\alpha_{ab}(z_{l_j})M_{ab}z_{l_j}\notag\\
&-\sum_{a=0}^{m-1}\sum_{b=0}^{m-a-1}\alpha_{ab}(z_{l_j})\begin{bmatrix}M_{a+1,b}\vspace{2mm} \\M_{a,b+1}\end{bmatrix}\notag\\
&+\sum_{a=0}^{m-1}\sum_{b=0}^{m-a-1}\beta_{ab}(z_{l_j})M_{ab}(v_{l_j}-\dot M_1)\big)-p\dot M_1\Big)
\label{HJB5}
\end{align}
with the boundary condition
\begin{equation}\tilde J(t_f,M_{ab}, \dot M_{10}, \dot M_{01}, z_l, v_l)=0.\end{equation}
\begin{remark}
If there are upper bounds for the errors we made during the approximations and if the original HJB equation (\ref{HJB3}) is stable backward in time, then the solution $\tilde J$ to (\ref{HJB5}) is close to the solution to (\ref{HJB3}).
\end{remark}

\subsection{Optimal control in an MPC setting}
Even though we can approximate the HJB equation so that it has all the moments and leaders' states as variables, solving this type of PDE numerically  is still difficult. There is no efficient method for solving general PDEs with more than four variables. At the same time, the approximation of the HJB equation performs poorly if the time horizon is big since the errors grow
larger in this case. We use model predictive control as a remedy. In particular, we optimize over a relatively short time period and apply the derived optimal control for only a few steps and repeat the procedure. At time $t=\tau$ and for a certain integer $q \in \mathbb{Z}$, we optimize over $q\Delta t$ time length in the future and apply the obtained optimal control for $\Delta t$. Then the sub-problem becomes:
\begin{equation}
\begin{array}{cl}\vspace{1mm}
\displaystyle\min_{u}~~ & \bar{\Phi}(z(\tau+q\Delta t))+\displaystyle\int_{\tau}^{\tau+q\Delta t}\bar{\mathcal L}(t)dt\\\vspace{1mm}
s.t.~~ &\dot z_i=v_i\\\vspace{1mm}
&\dot v_i=\sum_{j=1}^{M}\left(g_1(z_i,z_{l_j})(z_{l_j}-z_i)+g_2(z_i,z_{l_j})(v_{l_j}-v_i)\right)\\\vspace{1mm}
&~~~~~+f_i-p v_i,\text{ for }i=1,\cdots,N,\\\vspace{1mm}
 &\dot z_{l_j}=v_{l_j}\\\vspace{1mm}
&\dot v_{l_j}=u_j-pv_{l_j},\\\vspace{1mm}
&\|u_j\|\leq u_{max},~~~~~~~~~~\text{  for }j=1,\cdots,M,\\\vspace{1mm}
&z(\tau), v(\tau),z_l(\tau),v_l(\tau)\text{ given, }
\end{array}
\label{MPC1}
\end{equation}
Note that we have a terminal cost here and the cost function $\bar {\mathcal L}$ may also be different form the original problem. There is no universal way to model the new cost functions $\bar{\Phi}$ and $\bar {\mathcal L}$. Since we still want to avoid using individual state information for the followers, we should choose $\bar{\Phi}$ and $\bar {\mathcal L}$ to be functions of the moments only. Here, since the goal is still moments tracking, we keep using the old cost function $\bar{\mathcal L(t)}$ during the process, which equals to $\sum_{a+b\leq m}c_{ab}(M_{ab}(t)-M_{ab}^d(t))^2$. The terminal cost should be used to approximate the tail optimal cost of the states, which theoretically should be the solution of the HJB equation at time $t=t+q\Delta t$. However, this solution is unavailable beforehand. Here we add the final errors for the derivative of the first order moments, i.e., $(\dot M_{01}-\dot M_{01}^d)^2+(\dot M_{10}-\dot M_{10}^d)^2$, in addition to all the moments errors as an approximation of the tail cost.

The corresponding approximated HJB equation of problem (\ref{MPC1}) will be
\begin{align}
-\frac{\partial \tilde J}{\partial t}&= \bar{\mathcal L}(t) +\frac{\partial \tilde J}{\partial z_l}^T v_l-p \frac{\partial \tilde J}{\partial v_l}^T v_l-u_{max}\|\frac{\partial \tilde J}{\partial v_l}\|\notag\\
&+\sum_{a=0}^{m}\sum_{b=0}^{m-a}\frac{\partial \tilde{J}}{\partial M_{ab}}\begin{bmatrix}aM_{a-1,b} &bM_{a,b-1}\end{bmatrix}\dot M_1\notag\\
&+\begin{bmatrix}\frac{\partial \tilde{J}}{\partial \dot M_{10}}&\frac{\partial \tilde{J}}{\partial \dot M_{01}}\end{bmatrix} \Big(  \sum_{j=1}^{M} \big( \sum_{a=0}^{m-1}\sum_{b=0}^{m-a-1}\alpha_{ab}(z_{l_j})M_{ab}z_{l_j}\notag\\
&-\sum_{a=0}^{m-1}\sum_{b=0}^{m-a-1}\alpha_{ab}(z_{l_j})\begin{bmatrix}M_{a+1,b}\vspace{2mm} \\M_{a,b+1}\end{bmatrix}\notag\\
&+\sum_{a=0}^{m-1}\sum_{b=0}^{m-a-1}\beta_{ab}(z_{l_j})M_{ab}(v_{l_j}-\dot M_1)\big)-p\dot M_1\Big)
\label{HJBMPC}
\end{align}
with the boundary condition:
\begin{equation}
\tilde J(t+q\Delta t)=\bar{\Phi}(z).
\label{boundaryHJBMPC}
\end{equation}

The complexity of numerically solving the partial differential equation (\ref{HJBMPC}) grows exponentially with $q$, which means we cannot handle a long period prediction. On the other hand, the whole MPC approach may lead to an unstable closed-loop system if $q$ is too small. There is a tradeoff between the numerical complexity and the stability of the approach, which will be discussed and examined with a numerical experiment in the next section.

\section{Numerical experiments}

We will implement the optimal control approach on a specific model to show the capability and robustness of it. Let us use the following settings:
$$N=100,M=4,m=5,$$
which means there are 100 followers and 4 leaders. The moments we want to track is up to the fifth orders, i.e., $M_{ab}$ with $a+b\leq 5$ because we will deal with 4-polygon now. The reference signal for the moments are calculated from integrating $x^ay^b$ over a moving and shrinking square when assuming a uniform distribution. Namely, we want the followers uniformly spread inside a moving square while following the leaders. The center of the square should initially be located at the origin and slowly move to a point called ``exit" at the position of $\begin{bmatrix}120 & 120\end{bmatrix}^T$. In this case, there is a potential numerical problem of the approach in that for high order moments, the scale of their value becomes quite big. For example, the 5-order moment $M_{40}$ may become $120^{5}\approx 2\times 10^10$, which is much bigger than $M_{10}\approx 120$ when the followers reach the exit. Therefore, we change to a new coordinate system in each MPC iteration as below.

\subsection{Re-coordination and centralized moments}
When we want to solve (\ref{MPC1}) in each MPC iteration, we change the coordinates in order to make the current center-of-mass of the follower crowd, $\begin{bmatrix} M_{10} & M_{01}\end{bmatrix}^T$ to be the new origin. By doing this, the measured moments then become the centralized moments we defined in (\ref{moments2.2}). The reference signal needs to be updated as well and there are two ways doing that:
\begin{itemize}
{\item[1. ] Integrate the functions $x^ay^b$ over the original desired shape in the new coordinate system. Some extra calculations are needed in each iteration step.}
{\item[2. ] Track the original first order moments $M_{10}$ and $M_{01}$ in the new coordinate system while the other (centralized) moments track signals generated by integrating over a shape centered on the origin. Since those reference signals can be calculated off-line, this approach will speed up the calculation. However, the objective is changed from ``the crowd should follow `this' shape trajectory" to ``the center of the crowd should follow `this' path while the whole crowd should form `this' shape".}
\end{itemize}
The second approach alters the objective a little while not changing the ultimate goal of shape tracking, and is used in the simulation to simplify the computation.

We use the following functions for the model:
$$g_1(d)=0.3+20e^{-\frac{d}{10}}-\frac{20}{d+0.1}, \text{ and }g_2(d)=\frac{20}{d+0.1}.$$
$$f(z_i,v_i,z_k,v_k)=\frac{8}{d}e^{-0.2d}(z_i-z_k)$$
\begin{align*}
\bar{\Phi}=&c_0((\dot M_{10}^d-\dot M_{01})^2+(\dot M_{01}^d-\dot M_{01})^2)\notag\\
&~~~~~~~~~~~~+\sum_{a=0}^{m-1}\sum_{b=0}^{m-a-1}c_{ab}(M_{ab}^d- M_{ab})^2,
\end{align*}
$$\bar{\mathcal L(t)}=\sum_{a=0}^{m-1}\sum_{b=0}^{m-a-1}\bar c_{ab}(M_{ab}^d(t)- M_{ab})^2,$$
where $c_{ab}$ and $\bar c_{ab}$ can be tuned for the test purpose. Note that the moments we write here should already be transformed accordingly as we mentioned above.

We still need to do polynomial approximation of the functions $g_1$ and $g_2$ by using Taylor expansions. The expansion should be taken around the center-of-mass of the follower crowd, which is now the origin in the new coordinate system, and the degree of the polynomial should be $m-1=4$. The expansion will look like:
\begin{align*}
g_s(z,z_l)&\approx g_s(l)+g_{sx}(l)z_{x}+g_{sy}(l)z_{y}+\cdots \notag\\
&+\frac{1}{24}\Big(g_{sxxxx}(l)z_{x}^4+4g_{sxxxy}(l)z_{x}^3z_{y}\notag\\
&+6g_{sxxyy}(l)z_{x}^2z_{y}^2+4g_{sxxyy}(l)z_{x}z_{y}^3+g_{syyyy}(l)z_{y}^4\Big),
\end{align*}
where $l=\|z_l\|$, $g_{sx}=\frac{\partial g_s}{\partial x}$, $g_{sy}=\frac{\partial g_s}{\partial y}$, and the higher order partial derivatives are denoted in the same pattern, for $s=1,2$. Note that both the variables $z$, $z_l$ and the functions $g_1$, $g_2$ should be already transformed into the new coordinate system here.

\subsection{Numerical solution to HJB equation}
Once we have all the information needed, we should solve (\ref{HJBMPC}) with the boundary condition (\ref{boundaryHJBMPC}). We use the backward explicit finite difference method to solve the PDE by approximating the partial derivatives as:
$$\frac{\partial \tilde J(t+\Delta t)}{\partial t}\approx \frac{\tilde J(t+\Delta t)-\tilde J(t)}{\Delta t},$$
and
$$\frac{\partial \tilde J(\cdot,\xi)}{\partial \xi}\approx \frac{\tilde J(\cdot,\xi+\Delta \xi)-\tilde J(\cdot,\xi-\Delta \xi)}{2\Delta \xi},$$
where $\xi$ can be $M_{ab}$, $\dot M_{10}$, $\dot M_{01}$, $z_j$, or $v_j$. Then we can solve for the function value of $\tilde J$ at time $\tau$.

With the setting we made, there are 38 non-time variables for $\tilde J$ and each partial derivative needs function values at two points in the later time to approximate, which results to a complexity of $O(76^q)$ that grows exponentially with $q$ as mentioned above. Fortunately, it is still fast to solve with $q=3$. The simulations we make below assume $q=3$ and $\Delta t=0.1$.

\subsection{Simulation 1}

\begin{figure}
\centering
\subfigure[$t=0$]{
\label{Fig.sub.11}
\includegraphics[width=0.22\textwidth]{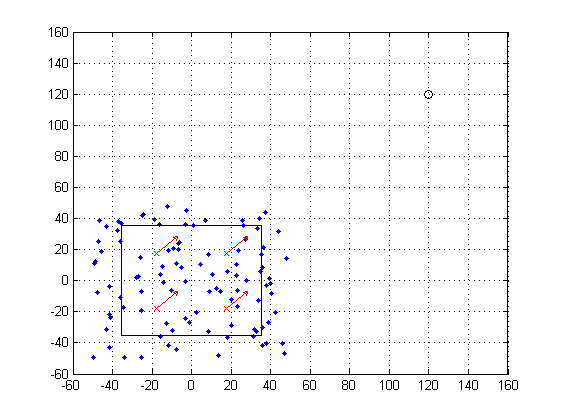}}
\subfigure[$t=20$]{
\label{Fig.sub.12}
\includegraphics[width=0.22\textwidth]{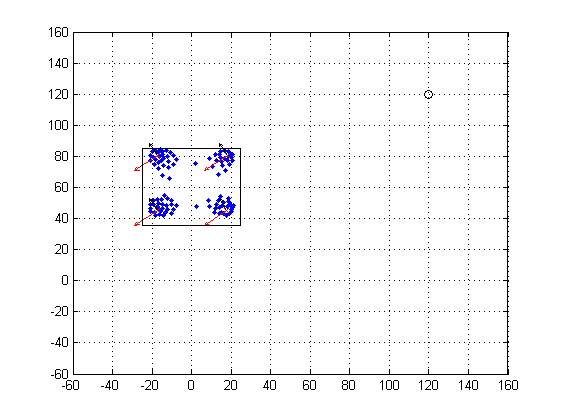}}
\subfigure[$t=40$]{
\label{Fig.sub.13}
\includegraphics[width=0.22\textwidth]{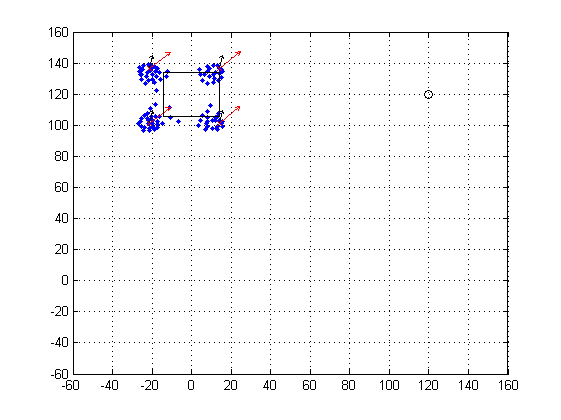}}
\subfigure[$t=80$]{
\label{Fig.sub.14}
\includegraphics[width=0.22\textwidth]{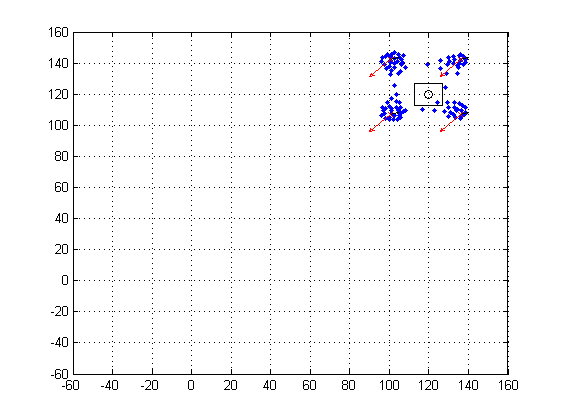}}
\caption{Snapshots for simulation 1. The blue dots indicate the position of the 100 followers and the red crosses are the 4 leaders. The squares are the desired shape in each plot. The black arrows are the velocity of the leaders while the red arrows are the acceleration of the leaders which come from the solution of (\ref{MPC1}).}
\label{sim1}
\end{figure}
In this simulation, the 100 followers are initially randomly distributed in the area of $[-50,50]\times[-50,50]$ while the 4 leaders start at $\begin{bmatrix}\pm25 & \pm25\end{bmatrix}^T$. The desired shape is a shrinking square first moving from the origin to the point $\begin{bmatrix}0 & 120\end{bmatrix}^T$, and then to the point $\begin{bmatrix}120 & 120\end{bmatrix}^T$ where the exit is. When the square reaches the exit, it will stay there with the size $10\times10$ for a while. We want the crowd to be uniformly spread in the desired shape. The coefficients in the cost functions are set to be \begin{equation} c_{ab}=\bar c_{ab}=\begin{cases}1000, &  k=1,\\10^{2(2-k)},&k>1,\end{cases}\end{equation}where $a+b=k$. The total length of the simulation is 120s.

\begin{figure}
\centering
\includegraphics[width=2.8in]{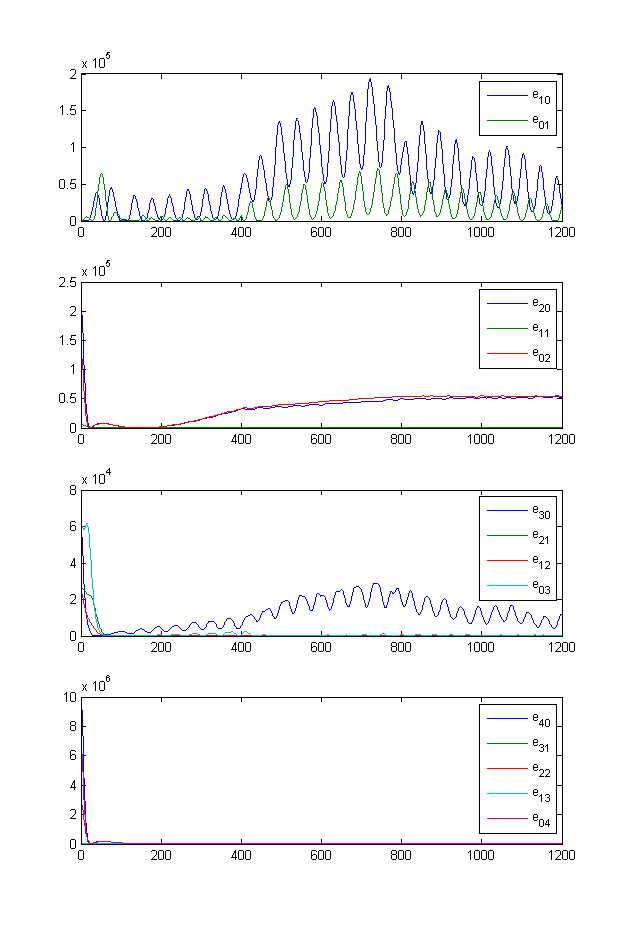}
\vspace{-10mm}\caption{Weighted tracking errors for different moments with the definition $e_{ab}=c_{ab}\|M_{ab}(t)-M^d_{ab}(t)\|^2$ up to the fourth order.}
\label{sim1e}
\end{figure}

Figure \ref{sim1} gives four snapshots of the simulation and Figure \ref{sim1e} gives the weighted tracking errors, where $e_{ab}=c_{ab}\|M_{ab}(t)-M^d_{ab}(t)\|^2$, up to the fourth order. The oscillation of the tracking error is expected since it is a double integrator model and the prediction horizon $q\Delta t$ is short.

\subsection{Simulation2}
In this experiment, we have the same setting of the followers and the leaders but add some obstacles in the environment as a disturbance to test the robustness of the method. The leaders will not take the extra disturbance from the obstacles into account when they calculate the optimal control, while both followers and leaders will react to the obstacles in their dynamics with an additional obstacle avoidance term:
\begin{equation*}
\dot v_i=\cdots+\sum_{k}h_k(z_i), ~~\dot v_{l_j}=\cdots+\sum_{k}h_k(z_{l_j}),
\end{equation*}
where $h_k$ models the obstacle avoidance force for the obstacle $k$. In this simulation, we set two round obstacles and model the forces as $$h_k(z)=\begin{cases}\frac{\kappa}{\|z-z_{ob_k}\|-r_k}(z-z_{ob_k}),&\text{if } \|z-z_{ob_k}\|\leq r_k+\delta,\\0,&\text{otherwise}\end{cases}$$where $z_{ob_k}$ is the center of the obstacle and $r_k$ is its radius, for $k=1,2$.  We choose
$$z_{ob_1}=\begin{bmatrix}5 \\ 60\end{bmatrix}, r_1=15, z_{ob_2}=\begin{bmatrix}80 \\ 80\end{bmatrix}, r_2=30,\delta=2.$$

\begin{figure}
\centering
\subfigure[$t=0$]{
\label{Fig.sub.21}
\includegraphics[width=0.22\textwidth]{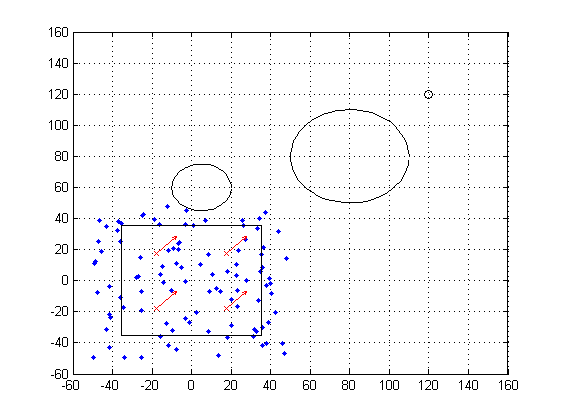}}
\subfigure[$t=12$]{
\label{Fig.sub.22}
\includegraphics[width=0.22\textwidth]{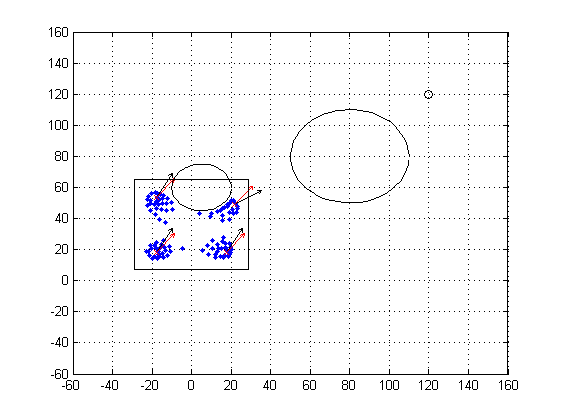}}
\subfigure[$t=20$]{
\label{Fig.sub.23}
\includegraphics[width=0.22\textwidth]{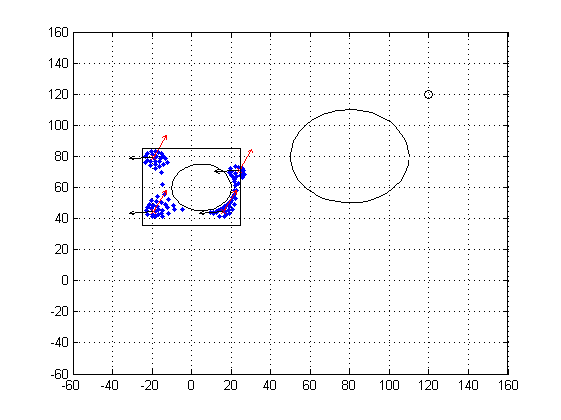}}
\subfigure[$t=30$]{
\label{Fig.sub.24}
\includegraphics[width=0.22\textwidth]{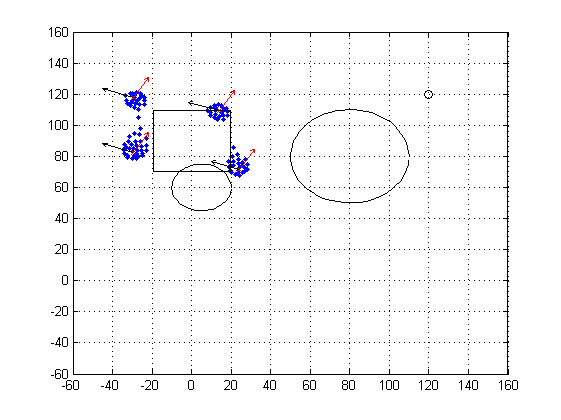}}
\subfigure[$t=60$]{
\label{Fig.sub.25}
\includegraphics[width=0.22\textwidth]{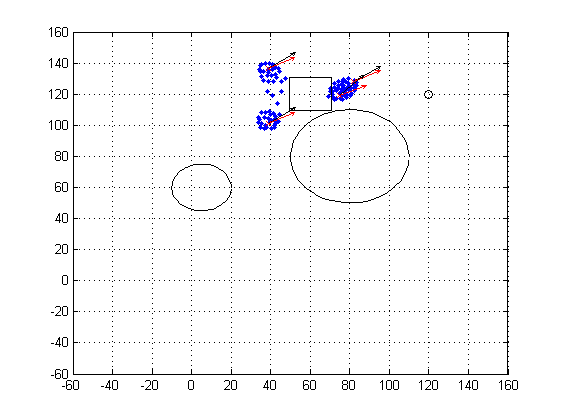}}
\subfigure[$t=120$]{
\label{Fig.sub.26}
\includegraphics[width=0.22\textwidth]{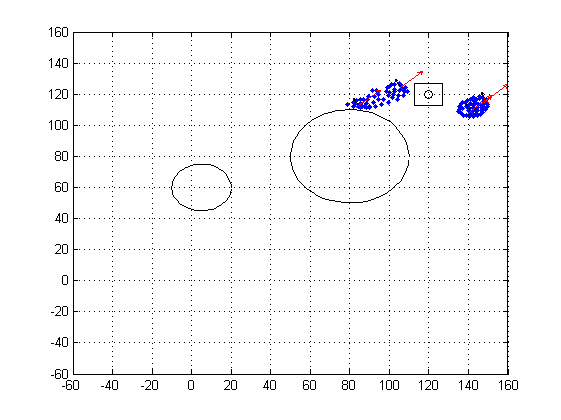}}
\caption{Snapshots for simulation 2. The blue dots indicate the position of the 100 followers and the red crosses are the 8 leaders. The squares are the desired shape in each plot. The black arrows are the velocity of the leaders while the red arrows are the acceleration of the leaders which come from the solution of (\ref{MPC1}). The two big circles illustrate the obstacles that both the followers and the leaders will avoid.}
\label{sim2}
\end{figure}
Figure \ref{sim2} shows that the method is still able to give a correct control signal when the followers and the leaders are passing by the obstacle.

\subsection{simulation 3}
\begin{figure}
\centering
\subfigure[$t=0$]{
\label{Fig.sub.21}
\includegraphics[width=0.22\textwidth]{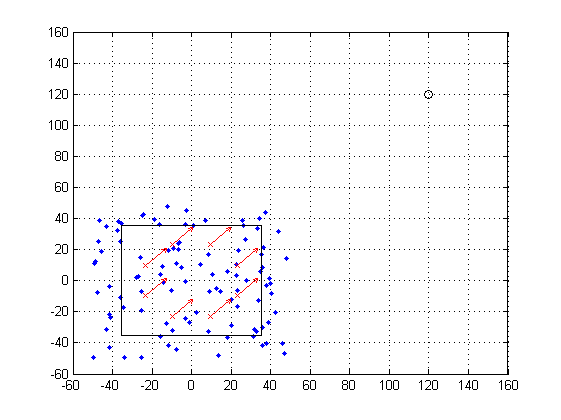}}
\subfigure[$t=20$]{
\label{Fig.sub.22}
\includegraphics[width=0.22\textwidth]{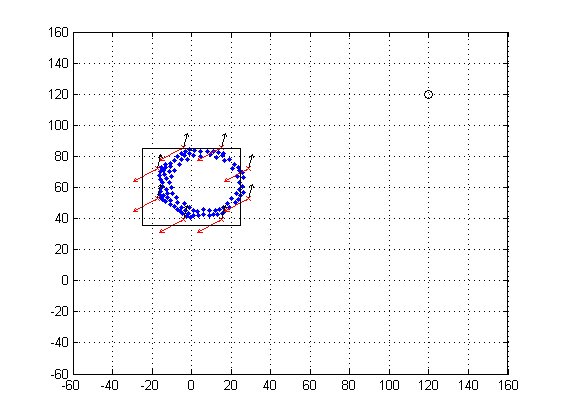}}
\subfigure[$t=40$]{
\label{Fig.sub.23}
\includegraphics[width=0.22\textwidth]{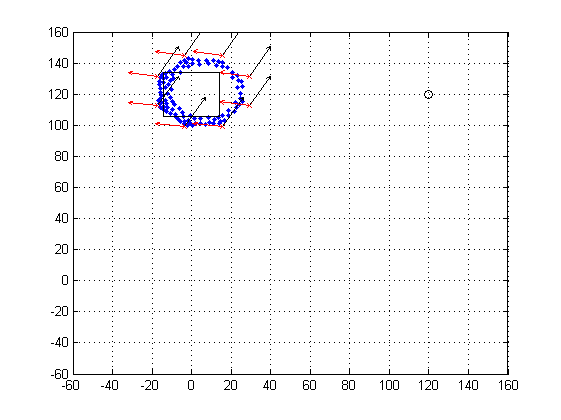}}
\subfigure[$t=120$]{
\label{Fig.sub.24}
\includegraphics[width=0.22\textwidth]{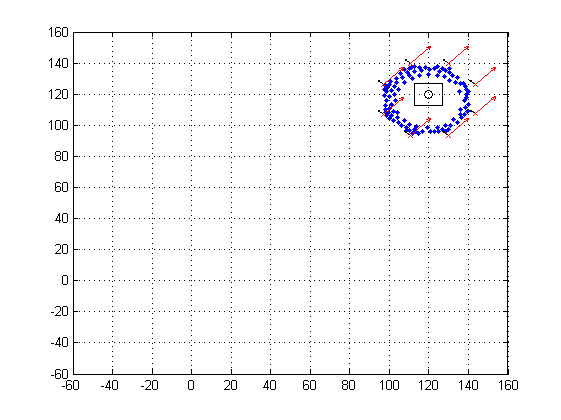}}
\caption{Snapshots for simulation 3. The blue dots indicate the position of the 100 followers and the red crosses are the 8 leaders. The squares are the desired shape in each plot. The black arrows are the velocity of the leaders while the red arrows are the acceleration of the leaders which come from the solution of (\ref{MPC1}).}
\label{sim3}
\end{figure}

\begin{figure}
\centering
\includegraphics[width=2.8in]{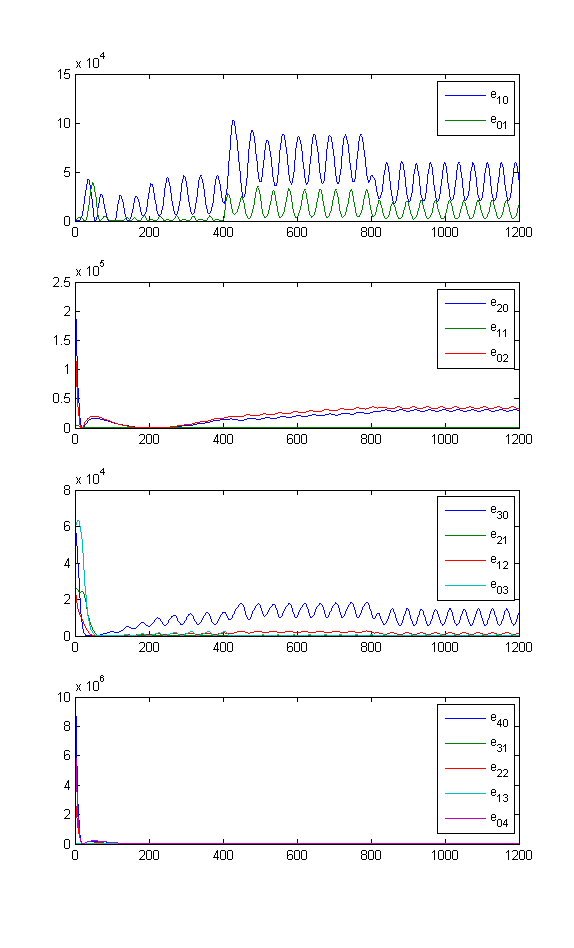}
\vspace{-10mm}\caption{Weighted tracking errors for different moments with the definition $e_{ab}=c_{ab}\|M_{ab}(t)-M^d_{ab}(t)\|^2$ up to the fourth order.}
\label{sim3e}
\end{figure}

In this experiment, we use almost the same setting as that in the first simulation except that we increase the number of leaders to 8.
The increase of the amount of leaders does decrease the errors slightly for all the orders. However, the performance improvement in terms of shape is hardly noticeable. In Figure \ref{sim3}, the followers form a circle in the end instead of the desired square, which could be considered to be worse than the shape generated by 4 leaders in Figure \ref{sim1} although the errors in Figure \ref{sim3e} turn to be smaller than those in Figure \ref{sim1e}.

\section{Conclusions and future work}

Since the moments of the crowd have a strong connection to its shape, we use moments to control the shape of the crowd. An optimal moments tracking problem is introduced in the paper and solved numerically with only using moments information as feedback. In order to further reduce the computational complexity, a model predictive control algorithm is used.  Three numerical experiments show that the method solves the moments tracking problem efficiently enough and can handle certain disturbance in the model. The performance of the optimal controller is acceptable even during the transient phase in the simulations.

Future work involves using local measurements of the moments to design a distributed leading strategy for the leaders. A better understanding of the relationship between the moments and shapes may give a more practical cost function to achieve a better performance in terms of shaping.

\end{document}